# Generalized Equal Area Criterion for Stability Analysis of Nonlinear Oscillators

Xin Xu, *Student Member, IEEE*, and Kai Sun, *Senior Member, IEEE*

*Abstract*— In power system analysis, the Equal Area Criterion is of great importance on revealing the condition for a grid-connected synchronous generator to maintain its transient rotor angle stability with the grid under a large disturbance. At present, the increasing renewable distributed energy resources (DERs) are introducing a variety of new nonlinear characteristics that are different from those of conventional generation although these DERs can still be operated and modeled as grid-connected oscillators with synthetic inertias. To study stability and synchronism of a grid-connected DER operated as a more general nonlinear oscillator, this paper proposes a Generalized Equal Area Criterion for estimation of transient stability margin subject to a large disturbance.

*Keywords*—Generalized Equal Area Criterion, nonlinear oscillator, transient stability.

## I. Introduction

Many real-life engineering systems can be modeled as networks of nonlinear oscillators and their normal operations and performances rely on synchronism and stability of the interconnected oscillators. For instance, for a power grid, the loss of synchronism among synchronous generators can result in rotor angle instability of generators, wide-spread electricity outages and economic losses across large areas [1]. Thus, the real-time stability analysis, monitoring and control on each oscillator of such systems subject to a major disturbance are of great importance and interests to both industry and academia. Compared to numerical simulation on all networked oscillators of a system, direct methods for transient stability analysis have their merits in fast time performance and quantitative margin information against instability, which are important for taking proper control actions before any detrimental consequences on the system.

In the field of power systems, a classical direct method for transient stability analysis on a synchronous generator connected to a power grid is the classical Equal Area Criterion (EAC), which treats the grid as an equivalent source or sink with a constant voltage so that the dynamics of the generator can be studied by a single-machine-infinite-bus (SMIB) system. The EAC tells that the generator will exit its domain of attraction to lose transient stability if its kinetic energy gained from a disturbance exceeds the potential energy it can obtain before passing an unstable equilibrium of the system [2]. The kinetic and potential energies can respectively be visualized as an "acceleration area" and a "deceleration area" respectively in a power-angle plane. The name "Equal-Area" comes from the fact that the system is certificated as stable if those two areas are equal. As a type of direct method by comparing an energy function (e.g. a Lyapunov function or its approximant) to its critical value, the EAC builds visual connection of the system's topological changes during and after a disturbance. Namely, the generator's rotor angle, mechanical power, electric power, energies and the equilibria of the system are all depicted in the same plane, making it easier to analyze transient stability by researchers and engineers. Hence, efforts were made to extend the EAC to analyze transient stability of multi-machine systems [3],[4].

However, the classical EAC has limitations when more renewables are integrated to the grid, especially the increasing power electronics interfaced distributed energy resources (DERs), which introduce new nonlinear characteristics different from those of conventional synchronous generators. These DERs can still be modeled as grid-connected oscillators, or even approximately as synchronous generators if they work in a synthetic inertia mode for frequency response purposes [5]. However, most of DERs introduce a variety of new nonlinear characteristics to power grids. For instance, a wind turbine is interfaced to the grid via a power electronic converter, which controls the power output in a more sophisticated manner than synchronous generators. Thus, the transient stability of a DER or non-conventional generator may not be analyzed using a classic EAC.

In this paper, a Generalized EAC (GEAC) is proposed for transient stability analysis of more general grid-connected nonlinear oscillators such as synchronous generators, DERs that operate as grid-connected oscillators, etc. It provides an approach for quantitative stability analysis and has the potential for on-line stability monitoring and control of, e.g. DERs in a power system. The rest of the paper is organized as follows. The methodology of the classical EAC is briefly introduced in section II for the comparison purpose. Then, the proposed GEAC and the associate approach for transient stability analysis are presented in section III. Case studies are conducted in section IV, and conclusions and future work are discussed in section V.

## II. Preliminary on Classical EAC

Consider an SMIB system that models a synchronous generator connected to a power grid. The generator can be considered as a nonlinear oscillator of one-degree-of-freedom (1-DOF), governed by swing equation (1):

$$\ddot{\delta} + \frac{D}{2H}\dot{\delta} + \frac{\omega_s}{2H}(P_{ek}(\delta) - P_m) = 0, \text{ where } P_{ek}(\delta) = P_{\max k}\sin\delta \quad (1)$$

This work was supported in NSF grants ECCS-1553863 and the ERC Program of the NSF and DOE under grant EEC-1041877.
X. Xu and K. Sun are with the Department of EECS, University of Tennessee, Knoxville, TN 37996 USA (e-mail: xxu30@vols.utk.edu, kaisun@utk.edu).

Rotor angle $\delta$ is the system state variable. $D$, $H$, $\omega_s$, $P_{maxk}$, $P_{ek}$ and $P_m$ are the damping coefficient, inertia, synchronous speed, maximum deliverable power, electrical power and mechanical power, respectively, and their units are second, second, rad/s, p.u., p.u. and p.u., respectively. Consider a fault near the generator and use subscript $k=1,2,3$ in (1) to denote the pre-fault, fault-on and post-fault conditions of the system, respectively. The fault only changes the maximum deliverable power $P_{maxk}$, while the other system parameters are unchanged.

Assume that the fault occurs at time $t = t_0$ and is cleared at time $t = t_c$. This scenario can be described as follows. Due to the fault occurrence and clearance, three periods of the system need to be considered along the time axis, namely, pre-fault ($t < t_0$), fault-on ($t_0 \leq t < t_c$), and post-fault ($t \geq t_c$) periods, which have different values of $P_{maxk}$.

To judge the stability of the generator, one need to study if rotor angle $\delta$ returns to the stable equilibrium point (SEP) of the post-fault system, say $\delta_{S,3}$. Without loss of generality, assume $P_{max2} < P_m < P_{max3} < P_{max1}$. Then, the stable equilibrium point (SEP) of (1), if exists, is $\delta = \delta_{S,k} < \pi/2$, and two unstable equilibrium points (UEPs) are $\delta_{U1,k} = \pi - \delta_{S,k}$ and $\delta_{U2,k} = -\pi - \delta_{S,k}$. Visualizing $P_m$ and $P_{ek}$ on the P-$\delta$ plane as shown in Fig. 1, we can see all the equilibria locate at the intersections of $P_m$ and $P_{ek}$. In terms of stability monitoring and assessment, one need to focus on if rotor angle $\delta$ returns to the SEP of the post-fault system, $\delta_{S,3}$.

The method of the classical EAC can be stated as follows. By ignoring the damping, i.e. $D=0$, kinetic energy $E_k$ is defined for (1) as the function of the power $P$ and rotor angle $\delta$, as the derivation shown in (2).

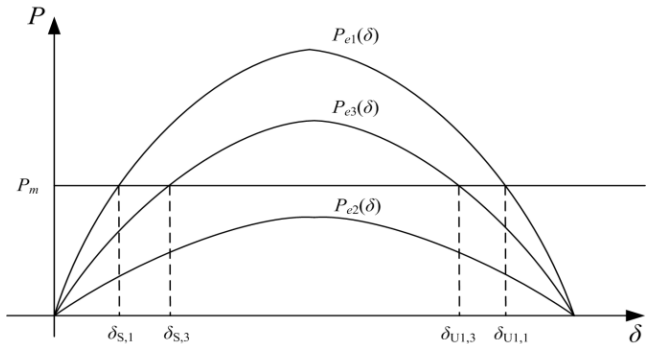

Figure 1.  $P_m$ and $P_{ek}$ on P-$\delta$ plane ($0 < \delta < \pi$)

$$\frac{2H}{\omega_s}\ddot{\delta}\dot{\delta} = [P_m - P_{ek}(\delta)]\dot{\delta}$$

$$\Leftrightarrow \frac{2H}{\omega_s}\frac{d}{dt}(\dot{\delta}^2) = [P_m - P_{ek}(\delta)]\dot{\delta} \quad (2)$$

$$\Leftrightarrow E_k = \frac{H}{\omega_s}\dot{\delta}^2 = \int[P_m - P_{ek}(\delta)]d\delta$$

Note that kinetic energy $E_k$ is equal to the "acceleration area" $A_{acc}$ on the P-$\delta$ plane as in Fig. 2, and the fault scenario can be visually described as follows. At $\delta(t = t_0) = \delta_{S,1}$, $P_{e1}(\delta_0)$ is decreased to $P_{e2}(\delta_0)$, and the kinetic energy is increased by area $A_{acc}$. The rotor is speeded up and the rotor angle $\delta$ begins to increase. At $\delta(t = t_c) = \delta_c$, $P_{e2}(\delta_c)$ is increased to $P_{e3}(\delta_c)$, the kinetic energy will decrease due to the "deceleration area" $A_{dec}$. The rotor will slow down before $\delta = \delta_{U1,3}$.

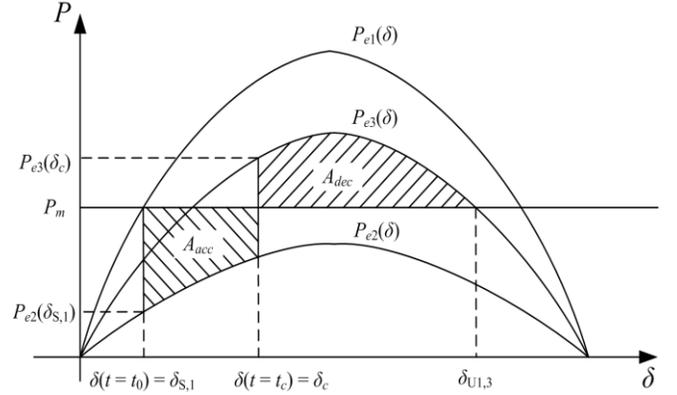

Figure 2.  Areas about $E_k$ on P-$\delta$ plane ($0 < \delta < \pi$)

Up to this end, let $\delta$ return to $\delta_{S,3}$. The EAC says $A_{dec} \geq A_{acc}$ in order for the system to maintain stability; otherwise, $\delta$ will pass $\delta_{U1,3}$, and the rotor will continue to accelerate and lose stability. Define this stability margin index:

$$Margin = \frac{A_{dec} - A_{acc}}{A_{acc}} \quad (3)$$

The kinetic energy defined in (2) and the corresponding "areas" in (3) are based on the fact that $P_{ek}(\delta)$ is a sinusoidal function with respect to $\delta$. However, this may unnecessarily be true on a DER with new nonlinear characteristics. Moreover, the assumption $D=0$ actually ignores the damping effect and makes the stability assessment result inaccurate. Therefore, there is a need to investigate a more generalized version of EAC in order to accommodate the new nonlinear characteristics and the damping effect.

III.  GENERALIZED EAC

A.  *System model*

Consider an autonomous 1-DOF nonlinear oscillator in (4) to generalize (1), whose SEP is $[\delta, \dot{\delta}]^T = [0,0]^T$.

$$\begin{cases} \dot{\delta} = \Delta\omega \\ \Delta\dot{\omega} = -a_0\Delta\omega - f(\delta) = -P_{\Delta\omega}(\Delta\omega) + P_f(\delta) \end{cases} \quad (4)$$

where $f(\delta) = a_{const}+a_1\delta+a_2\delta^2+a_3\delta^3+\ldots+ a_N\delta^N$, i.e. modeled by a truncated Taylor series. If $P_f > P_{\Delta\omega}$, $\Delta\omega$ will increase; otherwise it will decrease.

The basic idea behind the GEAC is stated as follows. A kinetic energy $E_g$ can be derived as in (5).

$$\Delta\dot{\omega}\Delta\omega = [-P_{\Delta\omega}(\Delta\omega) + P_f(\delta)]\Delta\omega$$

$$\Leftrightarrow \frac{1}{2}\cdot\frac{d}{dt}(\Delta\omega^2) = [-P_{\Delta\omega}(\Delta\omega) + P_f(\delta)]\Delta\omega \quad (5)$$

$$\Leftrightarrow E_g = \frac{1}{2}\Delta\omega^2 = \int[-P_{\Delta\omega}(\Delta\omega) + P_f(\delta)]d\delta$$

With the energy defined in (5), we can analyze the dynamics of the oscillator on the $P$–$\delta$ plane. On this plane, the curve of $P_f$ is fixed, while $P_{\Delta\omega}$ is initial value-dependent. First we will show how to qualitatively verify the system stability in the $P$–$\delta$ plane. Then, we define a quantitative criterion based on the energy in (5).

Note that on the $P$–$\delta$ plane, the rotor angle $\delta$ may reverse its moving direction, which is referred to as a turning point. The movement between two temporally adjacent turning points is referred to as a swing. A swing is considered "forward" if $\delta$ keeps increasing; otherwise it is a "backward" swing. Without loss of generality, only a "forward" swing is studied in the following sections.

### B. Qualitative Analysis

A qualitative analysis can be done by *(a)* monitoring $\delta$ and $P_{\Delta\omega}$, and *(b)* comparing $P_{\Delta\omega}$ with $P_f$ on the $P$–$\delta$ plane as in Fig. 3. The three roots of $P_f(\delta)=0$, namely, $\delta_1$, $\delta_2=0$, and $\delta_3$, need to be solved for comparison purpose. Note that $\delta_2$ corresponds to the SEP, while $\delta_1$ and $\delta_3$ correspond to the UEPs.

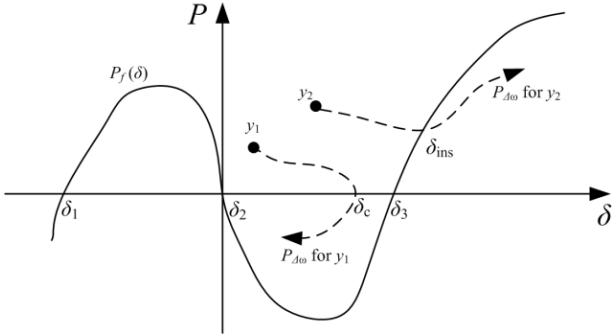

Figure 3.  Qualitative analysis on $P$-$\delta$ plane

In the "forward" swing, $\delta$ moves towards the right hand side. Two different points at the fault clearance moment are used to illustrate the qualitative analysis process for stable and unstable cases. One point is $y_1 = [\delta_{y1},\ P_\omega(\Delta\omega_{y1})]^T$, corresponding to the stable case, and the other one is $y_2 = [\delta_{y2}, P_{\Delta\omega}(\Delta\omega_{y2})]^T$, corresponding to the unstable case. The analysis is given below.

**Stable Case $y_1$:** As $\delta$ keeps increasing, both $\Delta\omega$ and $P_{\Delta\omega}$ keep decreasing due to $P_f < P_{\Delta\omega}$. At the moment that $P_{\Delta\omega}$ moves across the $\delta$-axis at $\delta=\delta_c<\delta_3$, both $\Delta\omega$ and $P_{\Delta\omega}$ turn negative, and $\delta$ begins to decrease. In this case, the rightmost UEP $\delta_3$ is not reached by the "forward" swing, and thus, that 'forward' swing is certificated to be stable.

**Unstable Case $y_2$:** As $\delta$ keeps increasing, both $\Delta\omega$ and $P_{\Delta\omega}$ keep decreasing due to $P_f < P_{\Delta\omega}$. At the moment that $P_{\Delta\omega}$ moves across $P_f$ at $\delta=\delta_{ins}$, $\Delta\omega$ begins to increase due to $P_f > P_{\Delta\omega}$, and $\delta$ also begins to increase towards positive infinity. In this case, by no means can the system return to the SEP, and the post-fault system loses stability in the "forward" swing.

For the unstable case, $\delta_{ins}$ can be used as an indicator of instability. Actually, the rightmost UEP $\delta_3$ verifies instability without the appearance of $\delta_{ins}$: if $P_{\Delta\omega}$ moves across the line $\delta=\delta_3$, it will definitely move across $P_f$ as in Fig. 3, and the system will lose stability in the "forward" swing.

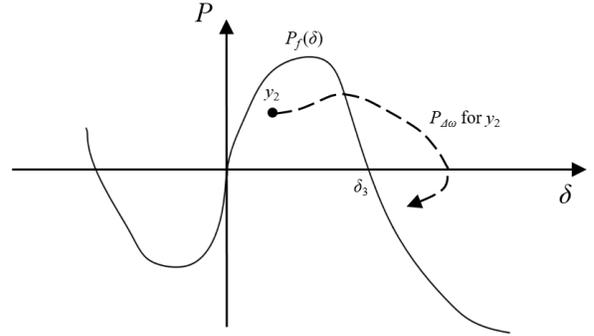

Figure 4.  Qualitative analysis: $P_f$ moves downwards at $\delta=\delta_3$

Note that in Fig. 3, $P_f$ moves upwards at $\delta=\delta_3$. However, if $P_f$ moves downwards at $\delta=\delta_3$ as it does in Fig. 4, $\Delta\omega$ and $P_{\Delta\omega}$ will eventually turn negative and $\delta$ will eventually reverse its moving direction. In other words, the system will never lose stability in the "forward" swing. Therefore, in this case, it is not necessary to monitor the stability for the "forward" swing.

In general, the qualitative analysis procedure to check the "forward" swing stability is stated as follows.

**Step 1**: Check if $P_f$ moves upwards or downwards when moving across the $\delta$-axis at the rightmost UEP $\delta_3$. If it moves downwards, the "forward" swing will never lose stability and we can start to assess the next swing; otherwise, go to step 2.

**Step 2**: Monitor the variation of $\delta$ and $P_{\Delta\omega}$. If $\delta$ becomes greater than $\delta_3$, the system loses stability; otherwise, if $P_{\Delta\omega}$ turns negative and $\delta$ reverses its moving direction, the system stays stable for the current swing, and we can start to assess the next swing.

The procedure for checking the "backward" swing stability can be similarly designed.

### C. Quantitative Stability Margin

Similar to the classical EAC, we can compare the "acceleration area" $A_{acc}$ and "deceleration area" $A_{dec}$ to estimate a quantitative stability margin. The methodology is illustrated by using the stable and unstable cases in Fig. 3.

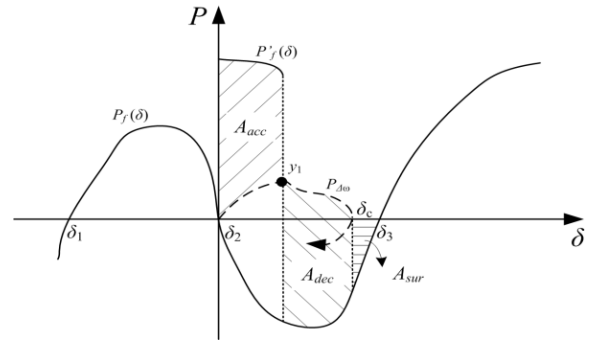

Figure 5.  Acceleration, deceleration and surplus area for stable case

**Stable Case $y_1$:** The "areas" for the stable case are illustrated in Fig. 5. Assume that, during the fault-on period, $P_f$ is temporarily switched to $P'_f$. The area between $P'_f$ and $P_{\Delta\omega}$ is $A_{acc}$, which accelerates $\delta$ to move rightwards. At the moment the fault is cleared, $P'_f$ is switched back to $P_f$. The area between $P_{\Delta\omega}$ and $P_f$ is $A_{dec}$. For stable case, it can be easily implied that $A_{acc} = A_{dec}$. Note that the area $A_{sur}$ marked in Fig. 5 is a redundant "deceleration" area and called the *surplus area*. With the areas defined above, the stability margin for the stable case is defined as:

$$Margin = \frac{A_{sur}}{A_{acc}} \quad (6)$$

**Unstable Case $y_2$:** The "areas" for the unstable case are illustrated in Fig. 6. Assume that, during the fault-on period, $P_f$ is temporarily switched to $P'_f$. The area between $P'_f$ and $P_{\Delta\omega}$ is $A_{acc}$, which accelerates $\delta$ to move rightwards. At the moment the fault is cleared, $P'_f$ is switched back to $P_f$. The area between $P_{\Delta\omega}$ and $P_f$ is $A_{dec}$. It can be easily implied that $A_{acc} > A_{dec}$. In this case the "deceleration area" is not adequate to let $\delta$ return. The stability margin is defined as:

$$Margin = \frac{A_{dec} - A_{acc}}{A_{acc}} \quad (7)$$

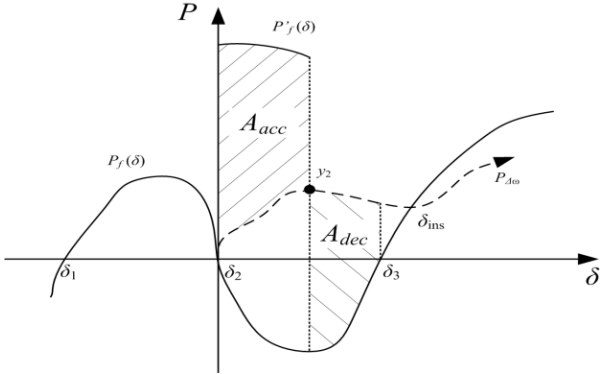

Figure 6. Acceleration and deceleration area for unstable case

For a "forward" swing, the stability margin is defined as:

$$Margin = \begin{cases} \dfrac{A_{sur}}{A_{acc}}, & \text{if } \delta \text{ reverses its direction} \\ \dfrac{A_{dec} - A_{acc}}{A_{acc}}, & \text{if } \delta > \delta_3 \end{cases} \quad (8)$$

### D. Relationship between Generalized and Classical EAC

The classical EAC can be viewed as a specific case of the GEAC: if we expand $P_{ek}(\delta)$ into a truncated Taylor series and redefine the SEP of (1) as the coordinate origin, (1) becomes exactly (4), so (1) is a specific case of (4). In addition, by taking $a_0=0$, it is easy to see that the kinetic energy definitions in (2) and (5) are actually consistent.

### E. Discussion on Online Application

If measurements of $\delta$ and $P_{\Delta\omega}$ are available online and the post-fault system model is known, the GEAC can be easily applied for a 1-DOF nonlinear system in the form of (4), like power systems. For each swing, the post-fault system can be analyzed via the procedure in section III.*A*, and the stability margin can be calculated by (8).

### *Remarks*

For many systems like power systems, the post-fault system model is usually known since the fault is cleared by certain protection measure, which results in a deterministic change in grid topology, e.g. line tripping. However, the fault-on system model is hard to be accurately obtained online. For instance, the short-circuit location on a transmission line is hard to be accurately located online. Hence, it is difficult to accurately obtain the $P'_f$ for the stability assessment of the first swing. This issue can be solved using the kinetic energy definition in (5): at the fault clearance moment $t = t_c$, the kinetic energy is totally determined by $\Delta\omega(t_c)$. Therefore, the $A_{acc}$ of the first swing can be computed by

$$A_{acc} = E_g(t_c) = \frac{1}{2}\Delta\omega^2(t_c) \quad (9)$$

## IV. CASE STUDIES

### A. 3rd order approximant of an SMIB system

The 3rd order approximant of a classical SMIB system is used as the post-fault system to test the proposed GEAC. By expanding the sinusoidal function of $P_{ek}(\delta)$ up to the 3rd order Taylor series truncation as in [6], the ODE expression with respect to the state variables $[\Delta\delta, \Delta\omega]^T$ is given in (10). Here $\Delta\delta$ is the rotor angle deviation in rad from the SEP. It has one SEP at $[0, 0]^T$ and two UEPs at $[1.9454, 0]^T$ and $[-3.0849, 0]^T$. The vector field in the phase plane is shown in Fig. 7(a).

$$\begin{cases} \dot\delta = \Delta\omega \\ \Delta\dot\omega = -4.42\times10^{-4}\Delta\omega - 0.2649\Delta\delta + 0.0503\Delta\delta^2 + 0.04414\Delta\delta^3 \end{cases} \quad (10)$$

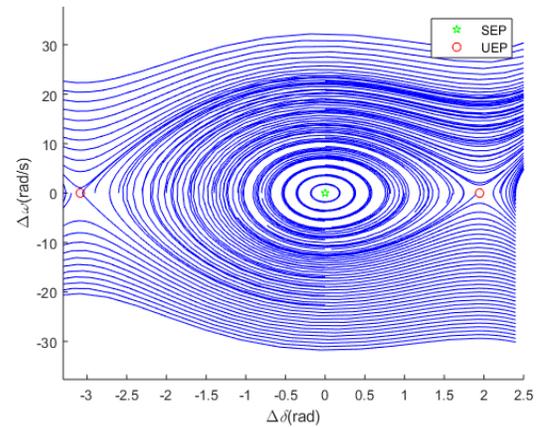

(a) Vector field in the phase plane

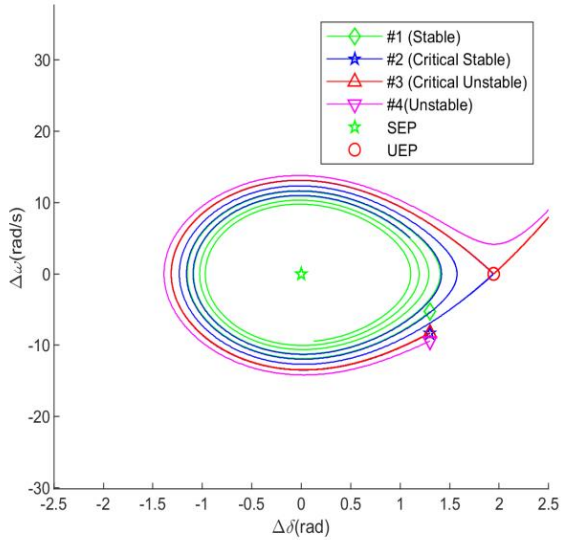

(b) system responses with different initial values

Figure 7. 3rd order approximant of an SMIB system

Use four different initial values at the fault clearance moment as shown in TABLE I, and assess the stability status of the post-fault system. The system responses starting from those initial values are shown in Fig. 7(b). We can see that as the initial value moves away from the SEP, the post-fault system tends to instability.

TABLE I. INITIAL VALUE AT FAULT CLEARANCE MOMENT

| #1 | #2 | #3 | #4 |
|---|---|---|---|
| $[0.13, -5.2779]^T$ | $[0.13, -8.3299]^T$ | $[0.13, -8.3315]^T$ | $[0.13, -9.4248]^T$ |

For each initial value, the stability assessment results for each swing are given in TABLE II, based on which the initial values are classified as stable, critical stable, critical unstable, and unstable cases. The "forward" swing is denoted by subscript F and "backward" by B. Note that the first swing is always a "backward" swing, i.e. $\delta$ decreases. The red mark means that the system loses stability at that swing.

The 2nd row in TABLE II contains the minimum stability margin of each case, clearly showing the trend toward instability as the initial value moves away from the SEP. The two critical cases have their stability margins close to zero. Thus, the proposed GEAC can accurately identify them.

In all cases, the "forward" swing is always more vulnerable to instability than the "backward" swing. Hence, it is suggested to take control measure for the "forward" swing to increase the stability margin, if possible. Meanwhile, if we only focus on the "forward" swings (or "backward") swings in the stable cases, we can see the stability margin keeps increasing as the time evolves. This is usually caused by the damping effect of the $\Delta\omega$ term.

### B. More general 3rd order oscillator

Consider a grid-connected DER to be operated as a nonlinear oscillator, which is not a good approximant of the SMIB system. We assume the coefficient of $\Delta\delta^2$ to change 0.0503 to 0.0603 as:

$$\begin{cases} \dot{\delta} = \Delta\omega \\ \Delta\dot{\omega} = -4.42\times10^{-4}\Delta\omega - 0.2649\Delta\delta + 0.0603\Delta\delta^2 + 0.04414\Delta\delta^3 \end{cases} \quad (11)$$

Use the same initial value in TABLE I to do the stability assessment. The system responses for each initial value are shown in Fig. 8. Note that the initial values #2-4 all lead to instability, which indicates the new devices might deteriorate the ability of the post-fault system to remain stable.

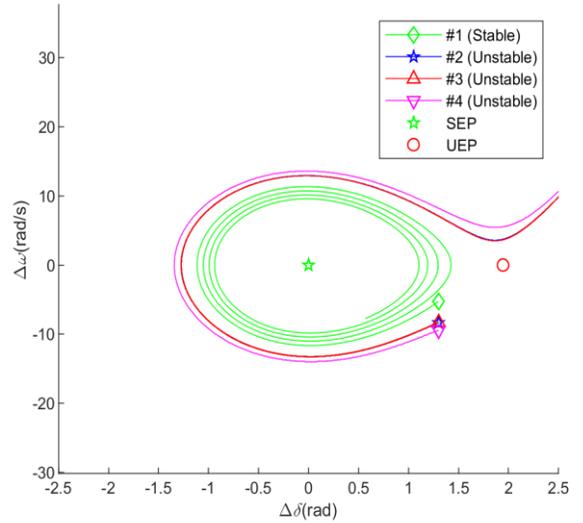

Figure 8. System responses of a more general 3rd order oscillator

For each initial value, the stability assessment results for different swings are given in TABLE III. Row two in TABLE III still contains the minimum stability margin of each case. The stability assessment result for the initial value #1 shows that the existence of the new devices generally decreases the stability margin of the "forward" swing but increases the stability margin of the 'backward' swing. The rest three cases go unstable since they do not have adequate stability margin in the first "forward" swing. Hence, extra control measure might be considered to improve the "forward" swing stability.

TABLE II. STABILITY ASSESSMENT RESULTS: SMIB

| No. of Swing | #1 (Stable) | #2 (Marginally Stable) | #3 (Marginally Unstable) | #4 (Unstable) |
|---|---|---|---|---|
| $1_B$ | 3.6701 | 3.3313 | 3.3311 | 3.1698 |
| $2_F$ | 0.2631 | 0.00005 | -0.0600 | -0.1499 |
| $3_B$ | 3.4372 | 2.4860 | - | - |
| $4_F$ | 0.4065 | 0.1192 | - | - |
| $5_B$ | 3.9569 | 2.9172 | - | - |
| $6_F$ | 0.5666 | 0.2490 | - | - |
| $7_B$ | 4.5369 | 3.3871 | - | - |
| $8_F$ | 0.7446 | 0.3924 | - | - |
| … | … | … | - | - |

TABLE III. STABILITY ASSESSMENT: SMIB WITH NEW DEVICES

| No. of Swing | #1 (Stable) | #2 (Unstable) | #3 (Unstable) | #4 (Unstable) |
|---|---|---|---|---|
| $1_B$ | 4.6069 | 4.2489 | 4.2534 | 4.0862 |
| $2_F$ | 0.1798 | -0.1351 | -0.1326 | -0.2172 |
| $3_B$ | 4.2444 | - | - | - |
| $4_F$ | 0.3142 | - | - | - |
| $5_B$ | 4.8655 | - | - | - |
| $6_F$ | 0.4643 | - | - | - |
| $7_B$ | 5.5558 | - | - | - |
| $8_F$ | 0.6305 | - | - | - |
| … | … | … | - | - |

## V. CONCLUSIONS AND FUTURE WORKS

A Generalized EAC is proposed for stability analysis of more general nonlinear 1-DOF oscillators and is compatible with the classical EAC. Future work will integrate this criterion into the online stability assessment of multi-oscillator systems by means of our recently proposed nonlinear model decoupling (NMD) methodology [7]-[10].